%% file: holext.tex
\documentclass[a4paper, 1pt]{article}
\usepackage{amssymb,amsbsy,amsmath,amsfonts,amssymb,amscd}
\usepackage{color,latexsym,eucal,exscale,epic,eepic,epsfig}
\date{}
\numberwithin{equation}{subsection}
\newtheorem{theorem}{Theorem}[subsection]
\newtheorem{lemma}{Lemma}[subsection]
\newtheorem{proposition}[theorem]{Proposition}
\newtheorem{corollary}[theorem]{Corollary}
\newtheorem{definition}{Definition\rm}[subsection]
\newtheorem{remark}[theorem]{Remark}

\def\ti{\widetilde }
\def\C{{\mathbb C}}
\def\NN{{\mathbb N}}
\def\M{{\mathcal M}}
\def\R{{\mathbb R}}
\def\V{{\mathcal V}}

\def\UU{{\mathfrak U}}
\def\CC{{\mathcal C}}
\def\n{{\mathfrak N}}
\def\m{{\mathfrak M}}

\def\l{{\ell }}
\def\B{{\mathcal B}}
\def\D{{\mathcal D}}
\def\O{{\mathcal O}}

\def\L{{\mathcal L}}
\def\U{{\mathcal U}}
\def\W{{\mathcal W}}

\begin{document}

\title{ON THE HOLOMORPHIC EXTENSION OF CR DISTRIBUTIONS FROM NON GENERIC
CR SUBMANIFOLDS OF $\C^L$}

\author {Nicolas Eisen~~eisen@math.univ-poitiers.fr}

\maketitle

\abstract{ \noindent We give a holomorphic extension result 
from non generic CR submanifold of $\C^L$ of positive CR dimension. We
consider $N$ a non generic CR submanifold given by $N=\{\n,h(\n)\}$
where
$\n$ is a generic submanifold of some $\C^{\ell}$ and $h$ is a CR map from
$\n$ into $\C^n$. We prove that if
$\n$ is a hypersurface then any CR distribution on $N$ extends holomorphically to 
a complex transversal
wedge, we then generalize this result for arbitrary $\n$ in the case
where the graphing function $h$ is decomposable at some $p'\in
\n$. We show that any CR distribution on $N$ that is
decomposable at $p=(p',h(p'))$  extends holomorphically to 
a complex transversal
wedge.}
\setcounter{section}{0}

\section{Introduction} 

\subsection {Statement of Results}
Let $\n$ be a smooth generic submanifold of $\C^{k+m}$ of CR dimension $k$
 and $h$ a smooth CR
map from $\n$ into some  $\C^{n}$ verifying $dh(0)=0$. Set $L=k+m+n$,
we construct a CR
submanifold $N$
of $\C^{p}$ near the origin as the graph of $h$ over $\n$, that is
$N=\{(\n,h(\n))\}$. It turns out that any non generic 
CR submanifold of $\C^L$
can be obtained in that fashion, see for example \cite{[Bo]}. The main 
question we address in this paper is the possible holomorphic
extension of CR functions of $N$ to some wedge $\W$ in a
 complex transverse direction.
We say that a vector in $\C^L$ is {\bf complex transversal} to $N$ at
$p \in N$ if $v \not \in $ span$_{\C}T_pN$.
For totally real submanifolds of $\C^L$, we have the following well known result (see the remarks following
 for different proofs
of this result). 
If $N$ is a non generic totally real 
submanifold of $\C^L$ and 
$ v \in \C^L$ is complex transversal to $N$ at $p$. Then for any
continuous function $f$ on a neighborhood of $p$ in $N$
 there exists $\W_v$, a wedge of direction
$v$  whose edge contains $N$
 such that $f$ has a holomorphic extension to $\W_v$.
In this paper, we shall study the case where $N$ is not totally
real, that is we consider CR submanifolds of $\C^L$ given by 
$N=\{(\n,h(\n))\}$   where $\n$ is generic in $\C^{k+m}$ of CR
dimension $k>0$. We say that a CR distribution extends holomorphically
to a wedge $\W$
if there exists a holomorphic function $F$ in $W$ such that the
boundary
value of $F$ on $N$ is $f$, i.e  $<f,\varphi>=
\lim_{\lambda \to 0^+} \int_{N}F(x+\lambda v)\varphi(x)dx$
 for any $v \in  \W$.
\begin{theorem} \label{th} Let $N=\{(\n,h(\n))\} $ 
be a smooth ($\CC^{\infty})$ non generic CR submanifold of
$\C^{k+1}\times \C^n$ 
such that $\n \subset \C^{k+1}$ is a hypersurface. If $f$ is a
CR distribution on $N$ then for any $p \in
N$ and any $v$ complex transversal to $N$ at $p$, there exists a wedge
$\W$ of direction $v$ whose edge contains a neighborhood of 
$p$ in $N$
and $F\in \O(\W)$ such that the boundary value of $F$ on $N$ is $f$.
\end{theorem}

We obtain some extension results for arbitrary $\n$ for decomposable CR distributions.
Since
the CR structure of $N$ is determined by $\n$, the CR distributions of
$N$ are precisely the CR distributions of $\n$.

\begin{definition}\label{D1} A  CR distribution $u$ on $\n$ is
decomposable at
the point $Z\in \n$ if, near $Z$, $u=\sum_{j=1}^KU_j$ where the $U_j$
are
CR distributions extending holomorphically to wedges $\W_j$ in $\C^{k+m}$  
with edges
$\n$. We shall say a  CR distribution  $u$ on $N$ 
is decomposable at a point $p=(p',h(p'))$ if $u$ is
decomposable at $p'$ on $\n$.
\end{definition}

\begin{theorem} \label{t1} Let $N=\{(\n,h(\n))\} $ 
be a smooth ($\CC^{\infty})$ non generic CR submanifold of $\C^{k+m+n}$ 
such that 
the function $h$ decomposable at some $p'_0 \in \n$.
Let  $v$ be a complex transversal vector to $N$ at $p_0=(p'_0,h(p'_0))$.

\medskip

{\bf (A)} If $f$ is a  decomposable CR distribution
 at $p_0$, then, near $p_0$, there exists a wedge
$\W$ of direction $v$ whose edge contains a neighborhood of 
$p_0$ in $N$
and $F\in \O(\W)$ such that the boundary value of $F$ on $N$ is $f$.

\medskip

{\bf (B)} There exists a wedge $\W$ of direction $v$
 and $\{F_{\l}\}_{\l=1}^n$, $F_{\l}\in
\O(\W)$ 
such that $dF_1\wedge ...\wedge dF_n\not = 0$ on $\W$ and each
$F_{\l}$ vanishes to order one on $N$.
\end{theorem}

We also get as a corollary:

\begin{corollary} \label{t3} Let $M$ be a smooth ($\CC^{\infty}$) 
generic submanifold of $\C^L$
containing through some $p_0 \in M$ a proper smooth ($\CC^{\infty}$) 
CR submanifold $N=(\n,h(\n))$ of
same CR dimension, write $p_0=(p'_0,h(p'_0))$ with $p'_0\in \n$. Assume that
the function $h$ decomposable at $p'_0$.
 Let 
$v\in T_{p_0}M \setminus 
\left [span_{\C}T_{p_0}N\right ]$. 
If $f$
is a decomposable CR distribution
at $p_0$, then there exists a wedge
$\W$ in $M$
of direction $v$  whose edge contains $N$
and $F$ smooth ($\CC^{\infty}$) CR functions on $\W$ such that
the boundary value of $F$ on $N$ is $f$.
 Furthermore, there exists a collection of smooth 
CR functions $\{g_{\l}\}_{\l=1}^n$
vanishing to  order one  on $N$ and 
such that $dg_1\wedge....\wedge
dg_n \not=0 $ on $\W$.
\end{corollary}

\subsection{Remarks} 
 
The equivalent of theorem \ref{th} for totally real manifolds
(mentioned
in the introduction)
 can be proved in several ways, one way is to follow
the theory of analytic vectors of an elliptic operator  due to  Treves 
(see section 3.2). Another way is to use the following result: Let
$N$ be a smooth submanifold of the boundary of $\Omega$, a strictly
pseudoconvex domain in $\C^L$. If $N$ is complex tangential ($TN
\subset (T(\partial \Omega) \cap iT(\partial \Omega))$) then $N$ is
a Pic interpolating set. See for example \cite{[Na]} or \cite{[Ru]}.
Given $N$, a totally real non generic submanifold of $\C^L$, one can
easily construct $\Omega$ as above and deduce the theorem.

\bigskip 

The boundary value of a holomorphic function is well defined for
functions of slow growth, that is
  there exists a constant $C>0$ and a positive integer
$\ell$
such that
$|F(z)|\leq {{C}\over{|dist(z,N)|^{\ell}}},$
where $dist(z,M)$ denotes the distance from a point $z$ to $N$.
 Conversely any CR distribution which extends holomorphically to a
 wedge $\W$ is given by the boundary value of a holomorphic function
 of slow growth. See
 for example \cite{[BER]}.

\bigskip 

As it is noted in \cite{[Tr-3]} on most CR submanifolds of $\C^L$
all CR functions are decomposable, hence the hypotheses of theorem
\ref{t1} hold in a generical sense. However, they are examples of
CR submanifolds of $\C^L$ on which indecomposable CR functions exist,
see \cite{[Tr-3]} for the general theory and \cite{[Rosay]} for an
elementary explanation of such examples. It turns out that on
Tr\'epreau's
original example of a CR manifold where undecomposable CR functions
exist
we still do get holomorphic extension of CR functions to complex
transversal
wedges (see section 3.2).

\bigskip

Theorem \ref{t1} implies that the extension obtained 
is not unique, which differs greatly with the holomorphic extension
results obtained for generic submanifolds.
We wish to point out the differences between our results and the
previously known results on holomorphic extension of CR functions.
The most general result on holomorphic extension is
Tumanov's theorem which states that if $M$ is a generic submanifold of
$\C^L$ that  is minimal at some $p_0$ then there exists a wedge $\W$
with edge $M$ (near $p_0$),
 such that one gets a ``forced'' unique local holomorphic
extension of CR functions defined in a neighborhood of $p_0$ in $M$.

\medskip

One should note that the question of CR extension can be viewed as a
Cauchy problem with Cauchy data on a characteristic set $N$. In
\cite{[Ei]}
we constructed an example of an abstract CR structure where there is
no such CR extension propriety. It is very easy to construct an 
abstract CR structure
where there is no CR function vanishing on $N$ to finite order
(see the example at the end of this paper).

\bigskip

Off course the situation is greatly different if we impose
holomorphic extension to a full neighborhood of $N$, as the 
next example shows.

\bigskip

\noindent {\bf Example} Let $I_n~~n>0$ be a sequence of disjointed
intervals (separated by some open sets) in $\R$ accumulating to the 
origin and let $u$ be a smooth
function such that $u|_{I_n}=1/n$. Let $\gamma=\{(s,u(s))\} \subset
\C^2$. Suppose $f$ is holomorphic on a neighborhood of the origin and
that $f|_{\gamma}=0$, then $f(w_1,{{1} \over {n}})=0$ for all $n$ for
$|w_1|$ small enough and thus $f \equiv 0$.

\subsection{Background}

The first result on holomorphic extension of CR function is due to
Lewy \cite{[Le]}. He proved that if a hypersurface was Levi 
nondegenerate at $p_0$
then CR functions extend holomorphically to one side of the
hypersurface. This result was generalized by Boggess and Polking
\cite{[Bo-Po]} for arbitrary dimensions. In the 
case of Levi flatness Tr\'epreau \cite{[Tr-1]}
proved that if a hypersurface in $\C^n$ is minimal at $p_0$ (it
contains through $p_0$ no $n-1$ dimensional complex manifold), then CR
functions extend holomorphically to one side of the hypersurface.
The generalization of Tr\'epreau's result to arbitrary codimension is
 due to Tumanov \cite{[Tu-1]}, in which he states that if the manifold $M$ is
minimal
at $p_0$ (it contains through $p_0$ no proper submanifolds of same
CR dimension) then CR functions extend to a wedge in ${\C}^n$ with
edge
$M$.
CR extension to manifolds of higher CR dimension has been dealt with
in works by Taiani and Hill \cite{[Hi-Ta]}, Tumanov
\cite{[Tu-2]}
as a well as in a recent paper by Eastwood and Graham \cite {[Ea-Gr]}.
 For
general background on CR geometry, we recommend the books by Baouendi,
Ebenfelt, Rothschild \cite{[BER]},
Boggess \cite{[Bo]} and Jacobowitz \cite{[Ja]} and the
 survey paper on holomorphic
extension by Tr\'epreau \cite{[Tr-2]}.

\subsection{Outline of the Paper}

The main part of this paper is the proof of theorem\ref{t1}, we then
obtain theorem \ref{th} as a corollary. We show that it suffices
to consider the case where $N$ is given as a trivial CR graph, that is
$N=\{(\n,0)\}$. We then solve a Dirichlet problem on some open subset
$\Omega$
of $\n \times \R^n$, by Treves's analytic vector theory, this yields
holomorphy
in the $w$ variables. To obtain holomorphy in $z$ we need to proceed
with
a deformation on $\n$ so that the ``boundary'' of $\Omega$ (boundary 
with respect to the variables of the elliptic operator) $\partial
\Omega$ is a CR submanifold of $\n \times \R^n$ of same CR dimension. 

\subsection{Acknowledgements}
The author would like to thank Jean-Pierre Rosay for some helpful
remarks and fruitful discussions as well as Jean-Michel Bony for the
references for the Dirichlet problem for
elliptic operators and Joel Merker for the
help with the pictures.


\section {Complex Transversal Holomorphic Extension in Manifolds of
type $N=(\n,0)$}

We first consider the case where the graph of $\n$ is trivial, i.e.  $h
\equiv 0$.
\begin{theorem}\label{Th2} Let $N=\{(\n,0)\} $ 
be a smooth ($\CC^{\infty})$ non generic CR submanifold of $\C^{k+m+n}$. 
Let  $\{v\}$ be complex transversal to $N$ at $p_0$. 
 If $f$ is a CR distribution
on $N$ that is decomposable at $p_0$, then there exists a wedge 
$\W$ of direction $v$ whose edge contains $N$
and $F\in \O(\W_{v})$ such that the boundary value of $F$ on $N$ is $f$.
\end{theorem}

The main tool to prove theorem \ref{Th2} is the following
proposition:

\begin{proposition}\label{Prop1} Let $N=\{(\n,0)\} $ 
be a smooth ($\CC^{\infty})$ non generic CR submanifold of
$\C^{k+m+n}$ and
let $f$ be a  distribution  on $N$ that admits a holomorphic
extension
to a wedge $\W$ in $\C^{k+m}$ with edge $\n$ near $p_0$,
then for any $v$ complex transversal to $N$ at $p_0$, there exists a
wedge
$\W_v$ in $\C^{k+m+n}$ of direction $v$ whose edge contains $\n$ and
$F \in \O(\W_v)$ such that the boundary value of $F$ on $N$ is $f$.
\end{proposition}

\noindent {\bf Proof of Theorem \ref{Th2}.} We assume proposition
\ref{Prop1}. Let $v$ be a complex transversal vector and let $u$ be a
CR
 distribution on $\n$. By hypothesis,
$u=\sum_{j=1}^{K}
U_j$ where each $U_j$ is a CR distribution that
 extends holomorphically to a wedge $\W_j\subset
\C^{k+m}$ with edge $\n$. We thus apply  proposition \ref{Prop1} to
each
$U_j$ to obtain a holomorphic extension to  wedges $\W'_j$ all in the
direction
$v$. Let $\W=\cap_{j=1}^K \W'_j$ we conclude that
the function $\sum_{j=1}^{K}
U_j$ is holomorphic in $\W$ and $\sum_{j=1}^{K}U_j=u$ on $\n$. 
This concludes the proof of theorem \ref{Th2}. $\blacksquare$

\bigskip

The remainder of this section will be devoted to the proof of proposition
\ref{Prop1}.

\subsection {Local Coordinates}

We begin with a choice of local coordinates on $\n$. 
$\n$ is a generic manifold in $\C^{k+m}$. We introduce local
coordinates near $p_0$. We may choose a local embedding so that
$p_0=0$ and 
$\n$ is parameterized in $\C^{k+m}=\C^k_z \times \C^m_{w'}$ by
 
\begin{equation} \label{e1} 
\n=\{(z,w')\in \C^k \times \C^m:Im(w')=a(z,
Re(w')),~~ a(0)=da(0)=0\}.
\end{equation}
We will denote by $s=Re(w') \in \R^m$, we thus have

\begin{equation} \label{e2} 
\n=\{(z,s+ia(z,s)\}\subset \C^k \times \C^m,~~~ T_0\n=\C^k\times \R^m.
\end{equation}

Define 
$\C T_p\n=T_p\n \otimes \C$ and
$T^{0,1}_p\n=T^{0,1}_p\C^{k+m} \cap \C T_p\n$. We say that $\n$ is a CR
manifold if dim$_{\C}T^{0,1}_p\n$ does not depend on $p$. The CR
vector fields of $\n$ are vector fields $L$ on $\n$ such that for any 
$p\in \n$ we have $L_p \in T^{0,1}_p\n$. One can
choose a basis $\L$ of $ T^{0,1}\n$ near the origin consisting of vector fields $L_j$
of the form

\begin{equation}\label{e3} 
L_j={{\partial} \over {\partial \overline z_j}}+\sum_{\ell=1}^nF_{j
\ell} {{\partial} \over {\partial  s_{\ell}}}.
\end{equation}

The wedge $\W$ in $\C^{k+m}$ with edge $\n$ on which we have
 holomorphic
extension of the function $f$  in a neighborhood of the origin is of the
form
$$
\W=\left (\U+i\Gamma\right ),
$$
where $\U$ is a neighborhood of the origin in $\n$ and $\Gamma$ is
 a conic neighborhood of some vector $\mu$  in $\R^m\setminus
 \{0\}$. 
We make another linear change of variables, identifying
 the vector $\mu$ with $(1,1,...,1)$ thus not changing \ref{e2}.
We therefore have for some $\eta>0$ and $B_{\eta}(0)$, a ball centered at the
 origin
 in $\C^{k+m}$ of radius $\eta$ 

\begin{equation} \label{e4}
 \begin{array}{l}
\W=\left (\U+i\Gamma\right ){\rm~~  ,where~} \U=B_{\eta}(0)\cap  \n~~
 {\rm and~} \Gamma~{\rm is~ a~conic~
 neighborhood~ of~}\\
 (1,1,....,1) ~{\rm in~} \R^m\setminus \{0\},\\
f~{\rm~extends~holomorphically~to~}\W=
 \left (\U+i\Gamma\right )\cap B_{\eta}(0).
\end{array}
\end{equation}

\subsection {Deformations}

 Denote by
$B_{r}$ the ball of radius $r$ centered at the origin in 
$\R^m$ and $\B_r$ the unit ball centered at the origin  in 
$\C^k$. Let $d \in \R$ be such that for $\eta$ given by \ref{e4} we have 
\begin{equation} \label{e7}
0<d<{{\eta}\over {\sqrt m}}.
\end{equation}
Let $\epsilon>0$ be given, let $\{b_j\}_{j=1}^m$
 be $\CC^{\infty}$ functions so that

\begin{equation} \label{e8}
\begin{cases}
b_j|_{\B_{\epsilon} \times B_{\epsilon}}=a_j|_{\B_{\epsilon} \times
B_{\epsilon}},\\
b_j=d~{\rm if}~s \in \R^{m} \setminus
B_{2\epsilon}
\end{cases}
\end{equation}

Define the generic submanifold of $\C^k \times \C^m$ $\ti \n$ as
follows

\begin{equation} \label{e9}
\ti \n=\{(z,s+ib(z,s)):(z,s) \in \C^k \times \R^m\}.
\end{equation}

\bigskip

\noindent {\bf Notation.} We shall use the following convention, if
$M$ is a manifold defined in a neighborhood of the origin,
parametrized
by $\C^k \times \R^{m}$, then  $M|_{B_r}$ is defined by

$$
M|_{B_r} =\{(Z(z,x)):x \in B_r\}.
$$
With the above notation we then have 
$$\ti \n|_{ B_{\epsilon}}= \n|_{ B_{\epsilon}},$$

which we illustrate with the following picture.

\bigskip

\begin{center}
\input p2.pstex_t
\end{center}

\begin{proposition}\label{p1} For $\epsilon$ small enough, there
exists
$b_j$ as in \ref{e8} so that $f$ (as in proposition \ref{Prop1}) extends
to a CR distribution on
$\ti \n$ in a $\eta$ neighborhood of the origin. 
\end{proposition}

\begin{corollary} \label{coro1}
 There is no loss of generality in assuming that the
functions $a$ in \ref{e1} verify \ref{e8}, i.e., one can replace $\n$
by $\ti \n$.
\end{corollary}

\noindent {\bf Proof of Proposition \ref{p1}.} 
The main tool is the following lemma.

\begin{lemma}\label{l1} There exists $b_j $ so that for
$\epsilon$ small enough  we have

$$\left (\ti \n \cap V_{\eta}\right ) \subset \overline \W,$$
where $\W$ and $V_{\eta}$ are as in \ref{e4}.
\end{lemma}
\noindent {\bf Proof of Lemma \ref{l1}.}

By \ref{e2} we have
\begin{equation} \label{e10}
\|a_j\|_{L^{\infty}(\B_{3\epsilon} \times B_{3\epsilon})}<C\epsilon^2.
\end{equation}

\bigskip

Choose $\vartheta$ and $\xi$ positive $\CC^{\infty}$ functions so that
 
\begin{equation} \label{e11}
\begin{cases}
\vartheta=\vartheta(s),~~~\xi=\xi(s),\\
\vartheta=0~{\rm on}~B_{\epsilon}~{\rm and}~\vartheta=1
~{\rm on}~\R^m \setminus B_{{{4}\over {3}} \epsilon},\\
\xi=1 ~{\rm on}~B_{{{4}\over {3}} \epsilon},
\xi=0 ~{\rm on}~~\R^m \setminus B_{{{5}\over {3}} \epsilon}.
\end{cases}
\end{equation}

Define the functions $b_j$ as follows,
$$b_1=d \vartheta+(1-\vartheta)a_1$$ and $b_j$ for $1<j\leq m$ is given
by
$$b_j=a_j+\xi\vartheta (d-a_1)+(1-\xi)(d-a_j),$$
where $d$ is given by \ref{e7}.
We see that
 $b_j=a_j$ if $s\in B_{\epsilon}$ and $b_j=d$ if $s \in B_{3\epsilon}
\setminus B_{2\epsilon}$.

We need to show now that a point $(z,s+ib(z,s))$ is in the wedge
$\W$. Write $ (z,s+ib(z,s))=(z,s+ia)+(0,iv)$ where $v$ is given
by
$$v=\vartheta(d-a_1) \left(1,...,1\right )+
(1-\xi)\left (0,a_1-a_2,...,a_1-a_j,...,a_1-a_m\right )$$

 By \ref{e10} $a_j=O(\epsilon^2)$, so if $\epsilon$ is small enough
we see that $ (z,s+ib(z,s))\in \W$
 To conclude the
proof
of lemma \ref{l1}, we need to make sure that $\epsilon$ is also small
enough so that $(z,s+id) \in B_{\eta}(0)
$ when $(z,s)\in \B_{3\epsilon}
\times B_{3\epsilon}$,
 which we can do by \ref{e7}.
 $\blacksquare$

\bigskip

The proof of  proposition \ref{p1} is now immediate, the CR extension
  of $f$ to $\ti \n$ is given by the restriction of the holomorphic extension of $f$
 to  $\ti
\n$, which is possible by the preceding lemma. $\blacksquare$

\bigskip

\noindent {\bf (Abuse of) Notation.} From henceforth we shall be working on
$\ti \n$. However, to keep the notation as simple as
possible
we will forget
the tilde and keep the $\n$.

\subsection{Resolution of a Dirichlet Problem}

 Consider the generic manifold $\n \times
\R^n \subset \C^{k+m}\times \C^n$.

\begin{lemma}\label{l2}
A basis $\L$ of the CR vector fields $\n \times
\R^n$
   consists of vector field $ L_j$ of the form

\begin{equation}\label{e17}
  L_j={{\partial} \over {\partial \overline z_j}}+\sum_{l=1}^mf_{j,l}
{{\partial} \over {\partial s_l}},
\end{equation}

where the  $  L_j$ can be chosen to verify the following property.

\medskip

If $s \in \R^{m} \setminus B_{2\epsilon}$ we have
$$  L_j={{\partial} \over {\partial \overline z_j}}.$$

\end{lemma}

The proof of the first part of lemma \ref{l2} can be found in any
textbook
on CR geometry. The second part is just a consequence
of the construction of $  \n$.

\subsubsection{Construction of the operator $\Delta_{\n\times\R^n}$}

\begin{lemma}\label{l3} There exists $m$ vector fields $R_j$ of the
form
$$R_j=\sum_{l=1}^m a_{jl}(z,s){{\partial} \over {\partial  s_l}}$$
where the $a_{jl}$  are smooth functions, such that

\medskip
 
(i) $R_{j}(  w'_l)=\delta_{jl}$, $j,l \in \{1,...,m\}$,

\medskip

(ii) $[R_j,R_l]=0$,
 
\medskip

(iii) $[  L_j,R_l]=0$,

\medskip

(iv) the set $\{  L_1,...,  L_k, \overline{   L_1},...,\overline
{  L_k},R_1,...,R_{m}\}$ span the complex tangent plane to 
$  \n $ near the origin.
\end{lemma}

The proof of the lemma is classic. (ii) (iii) and (iv) are a consequence of (i). We thus
 determine the $R_j$ by solving
for their coefficients in (i) (see for example \cite{[BER]} Lemma
8.7.13 page 234).

\begin{remark} \label{r1} We note that the $R_j$'s verifies the 
following property,

if $j \leq m$ then 
$$R_j|_{0}= R_j|_{\C^k\times \left [\R^{m}\setminus
B_{2\epsilon}\right ]}={{\partial}\over {\partial s_j}}. $$
 
\end{remark}

The equivalent of lemma \ref{l3} for $\n\times\R^n$, is given by

\begin{lemma} \label{l'3} The $m+n$ vector fields given by
$\{R_1,....,R_m,{{\partial}\over{\partial
t_1}},...,{{\partial}\over{\partial t_n}}\}$ verify the same
properties as in lemma \ref{l3}, where we replace $  \n$ by
$\n\times\R^n$.
\end{lemma}

 Define $\Delta_{\n\times\R^n}$ by

\begin{equation}\label{e1.1}
\Delta_{\n\times\R^n}=\sum_{j=1}^{m}R_j^2 +\sum_{j=1}^{n}
{{\partial^2}\over{{\partial
t_j}}^2} .
\end{equation}

From \ref{r1} and lemma \ref{l3} (iii) we deduce immediately the following:

\begin{lemma} \label{p2} The operator $-\Delta_{\n\times\R^n}$ is strongly
elliptic of degree two on $\n\times\R^n$  with
smooth
($\CC^{\infty}$) coefficients and no constant terms.
\end{lemma}

$\Delta_{\n\times\R^n}$ is a differential operator in the variables $s$ and
$t$
whose coefficient functions depend smoothly on the $z$ variable. We
are going to study a Dirichlet problem on $\n\times\R^n$. To do this, we
need 
$\Omega$ some open set in $\n\times\R^n$ with boundary $\partial \Omega$
parametrized
by some closed submanifold of $\R^{n+m}$ of codimension one (recall
that $\Delta_{\n\times\R^n}$ is a differential operator in the variables $s$ and
$t$).

\begin{remark} 
The reason we went to the trouble of redefining $  \n$ was to be able
to
define an open set $\Omega$ whose  boundary is a CR submanifold of
$\n\times\R^n$
of same CR dimension containing near the origin $\n$.
\end{remark}

\subsubsection{Construction of the open set $\Omega$}

In $\R^{m+n}$ let $\omega$ be an open set contained in 
$\R^{m+n} \cap \{t_1>0\}$ and whose boundary (see figure) is the union
of
two submanifolds $\gamma_1$ and $\gamma_2$, verifying the following
properties.

\bigskip

\noindent {\bf(1)} $\gamma_1$ is given by $t_1=0$ 
for $(s,t)\in B_{2\epsilon}$.

\bigskip

\noindent {\bf (2)} 
$\gamma_2 \subset B_{3\epsilon}\setminus B_{2\epsilon}$ is such that 

\medskip

(i) $\gamma_1 \cup \gamma_2$ bounds an open set $\omega$,

\medskip

(ii) $\gamma_1 \cup \gamma_2$ is smooth.

\bigskip

Define $\Omega$ by 
$$
\Omega=\{(z,  w(z,s,t): z\in \B_{3\epsilon}, (s,t) \in \omega\}.
$$
The boundary of $\Omega$ on which we shall impose the Dirichlet data
is defined to be 
$$\partial \Omega=\{(z ,  w(z,s,t):z\in \B_{3\epsilon} ,(s,t) \in  \gamma_1 \cup \gamma_2
\}.
$$

\begin{lemma}\label{l4} $\partial \Omega$ is a smooth CR submanifold
of $\n \times \R^n $ of same CR dimension.
\end{lemma}
\noindent{\bf Proof of Lemma \ref{l4}.} 
The lemma is clear for the part of $\gamma_1$ which contains $\n
\times \R^{n-1}$. Out of $\n
\times \R^{n-1}$ by lemma\ref{l2}
  the
generators of the CR vector fields $L_j$  are equal
to ${{\partial} \over {\partial \overline z_j}}$, hence the lemma.
$\blacksquare$

\bigskip

We have the following picture.

\bigskip

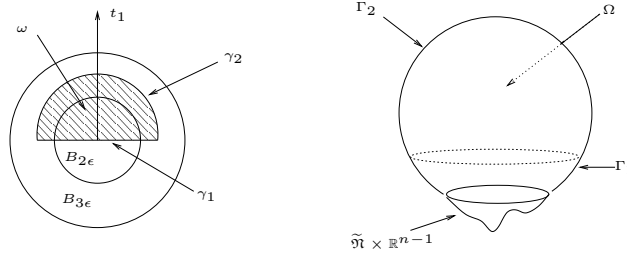
\begin{center}
\input p4.pstex_t
\end{center}

Denote by $(D)$ the Dirichlet problem on $\Omega$, i.e.
$${(D)}~~~~~~~~~~~~~~~~~
\begin{cases}
\Delta_{\n\times\R^n} (u)=0~~{\rm in}~\Omega,\\
u=g~~{\rm on}~\partial \Omega.
\end{cases}$$

We now have the following result, which seems well known but is not so
easy to find in the literature.
\begin{theorem}\label{proba} 
 For $g \in \D'(\partial \Omega)$, $(D)$ has a unique solution.
\end{theorem}
For a reference, we send the reader to \cite{[C-P]} theorem 5.2 and
remark 5.3 page 263-265, as well as \cite{[Ho]} chapter 10.

Let $S$ be the solution operator for $(D)$, that is for 
$g \in \D'(\partial \Omega)$, we have

\begin{equation}\label{ee2}
\begin{cases}
\begin{array}{l}
\Delta_{\n\times\R^n} (S(g))=0~~{\rm in}~~\Omega,\\
S(g)=g ~~{\rm on}~\partial \Omega.
\end{array}
\end{cases}
\end{equation}

\begin{proposition}\label{p3} Let $  L_j \in  
\L$ and let $g$ be a  CR distribution on $\partial {\Omega}$  then we
have
$$  L_j(S(g))=0.$$
\end{proposition}

\begin{lemma}\label{LL1} If $f$ is a CR distribution on  $\partial \Omega $
  then for any $  L_j \in   \L$ we
have 
$$
  L_j(S(f))|_{\partial \Omega}=0.
$$
\end{lemma}

\bigskip

\noindent {\bf Proof of Lemma \ref{LL1}.} The lemma holds since
$\partial \Omega$ is a CR submanifold of  $\n\times\R^n$
and $f$ is a CR distribution on $\partial \Omega$,  so 
$ L_j(S(f))|_{\partial \Omega}=  L|_{\partial \Omega}(f|_{\partial \Omega})=0$.

\bigskip

\noindent {\bf Proof of Proposition \ref{p3}.}
Consider $  L_j(S(g))$. By lemma
\ref{l3}
(iii) we have 
$$\Delta_{\n\times\R^n}   L_j(S(g))=  L_j\Delta_{\n\times\R^n}(
S(g)),$$
 hence we
have
\begin{equation}\label{e1.3}
\Delta_{\n\times\R^n} \left [  L_j(S(g))\right]=0.
\end{equation}
 
By lemma \ref{LL1} $  L_j(S(g))$ vanishes on $\partial \Omega$,
 hence by theorem \ref{proba} we have
 $  L_j(S(g))=0$ in $\Omega$. $\blacksquare$

\subsection{Analytic Vector Theory}

In this section we present the results developed by Baouendi and Treves
\cite{[Ba-Tr]} and later by Treves in \cite{[Tre]}. We have included this
section
for the sake of completeness and claim no originality whatsoever. We
will
state the main results needed here, sending the reader who wishes
to read the details of the proofs to  \cite{[Tre]}. 

\begin{definition} Using the notations
$R_{m+j}={{\partial}\over{\partial t_j}}$ and
$R^{\alpha}=R_1^{\alpha_1}...R_{m+n}^{\alpha_{m+n}}$ for $\alpha \in
\NN^{m+n}$,
we shall say that a continuous function $f$ in $\omega$ is an analytic
vector
of the system of vector fields $\{R_1,...,R_{n+m}\}$ if 
$R^{\alpha}f \in \CC^0 $ for any  $\alpha \in
\NN^{n+m}$, and if to every compact set $K$ of $\omega$ there is a
constant
$\rho>0$ such that, in $K$,
\begin{equation}\label{2.2.1}
\sup_{\alpha \in
\NN^{n+m}} \left ( \rho^{|\alpha|} {{ |R^{\alpha}f|} \over {\alpha!}}
\right )<\infty.
\end{equation}
\end{definition} 

Let $p
=(z,  w(z,s_0,t_0)) \in \Omega$ so that with our previous notation
$(s_0,t_0) \in \omega$. We are using the following convention

$$(z,  w(z,s_0,t_0)) \in \Omega \Leftrightarrow (z,s_0,t_0) \in
\B_{3\epsilon}\times \omega.$$

\begin{proposition}\label{p5} (Lemma 4.1 in \cite{[Ba-Tr]}). Let $p
=(z,  w(z,s_0,t_0)) \in \Omega$, let $\U$
be an open neighborhood of $(s_0,t_0)$ in $ \omega$, $K$ a compact 
subset of $\U$.
 To every pair of positive constants $C_1, C_2$, there
is
another pair of positive constants $C_1',C_2'$ with $C_2'$ independent  
of $C_1$ such that for any $\CC^{\infty}$ function $f$ in an open
neighborhood
of the closure of $\U$, if we have, for all $k\in \NN$ and all $z\in
\B_{3\epsilon}$,
\begin{equation}\label{e3.1}
\|\Delta_{\n\times\R^n} ^k f\|_{L^{\infty}(\U)} \leq C_1 C_2^k (2k)!,
\end{equation}
then for every $z \in \B_{3\epsilon}$ and every $\alpha \in \NN^{n+m}$
\begin{equation}\label{e3.2}
\|R^{\alpha} f\|_{L^{\infty}(K)}\leq C_1'{C_2'}^{|\alpha|}|\alpha|!
\end{equation}
\end{proposition}

\begin{proposition}\label{p4} (Proposition II.4.1 in \cite{[Tre]}). Let
$p \in \Omega$,  in
order that $f\in \CC^0$ be an analytic vector of the system of vector
fields $\{R_1,...,R_{n+m}\}$ it is necessary and sufficient that there
exists an open neighborhood $\V$ of $p$ in $\C^k \times \C^{n+m}$ and
a continuous function $
F(z,w)$ in $\V$ holomorphic with respect to $w$
and such that $f(z,s,t)=F(z,  w(z,s,t))$.
\end{proposition}

The main difficulty in the proof of proposition \ref{p4} is to show
that
the function defined by
$$F(z,  w,s,t)=\sum_{\alpha \in \NN^{m+n}}{{R^{\alpha}f(z,s,t)} \over
{\alpha!}}
\left (  w-  w(z,s,t) \right )^{\alpha}$$
is equal to $f$ for $  w$ near $  w(z,s,t)$ in $\Omega$ if $f$ is
an analytic vector of the vector
fields $\{R_1,...,R_{n+m}\}$.

\bigskip

We shall use proposition \ref{p5} to construct analytic vectors of the
vector
fields $\{R_1,...,R_{n+m}\}$  and then apply proposition \ref{p4}
to these vectors to obtain a holomorphic extension in the variables
$  w$.

\subsection{Proof of Proposition \ref{Prop1}}
 
We have
\begin{equation}\label{e'5} 
T_0N=\C^k\times \R^m \times \{0\}\subset \C^k\times \C^m\times\C^n.
\end{equation}
If $v$ is complex transversal to $N$ at the origin, then $v\in\{0\}
\times \{0\}  \times\C^n$. Therefore after a linear change of variables in 
$\C^{k+m+n}$ (viewed as a linear change of variables in $\C^n$ being
the identity on $\C^{k+m}$ so that \ref{e4} is unaffected) we may
assume
that $v$ has the following form
\begin{equation}\label{cpx} 
v=(0,0,v'),~~v'=(1,0,...,0) \in \R^n \times\{0\}\subset \C^n.
\end{equation}

We first show that $f$ extends holomorphically to a neighborhood of $\Omega$.

\bigskip

\noindent {\bf Notation} Given $P=(z,  w)$ a point in $\Omega$,
 $\U$ will denote a neighborhood of $  w$ in $\Omega$
and $U$ will denote a neighborhood of $  w$ in $\C^{m+n}$.

\bigskip

\noindent {\bf (A) Holomorphy in
  $  w$}

\bigskip

By construction, we have
$\Delta_{\n\times\R^n} S(f)=0$ on $\Omega$. Let 
$Q=(z_0,  w_0) \in \Omega$.
 By proposition \ref{p5}, $S(f)$ is an analytic vector for the
vector fields $R_j$ on some open neighborhood $\U_1$ of $  w_0$.
By proposition
\ref{p4},
$S(f)$ extends holomorphically as a function of $  w$ to a function
denoted by $  F_1$
 when $(z,  w)$ is in some open  neighborhood $\B_{3\epsilon}\times
U_1$
of $(z_0,  w_0)$.
  Let $  w_1 \in \U_1$, $  w_1 \not =  w_0$. Consider then 
$P=(z_0,  w_1)$.
By the above reasoning, we find the existence of neighborhoods
$\U_2$ and $U_2$
of $  w_1 $  and a function $  F_2$ defined on $\B_{3\epsilon} \times
U_2$, holomorphic as a function of $  w$ such that 
$$S(f)=F_2|_{\B_{3\epsilon} \times
\U_2}.$$

\bigskip

\noindent {\bf Claim:}   $ 
F_1=  F_2$ on $\B_{3\epsilon} \times \left [U_1 \cap U_2\right ]$.

\bigskip

\noindent {\bf Proof of the Claim.} 
For $z$ fixed,  $(\U_1 \cap \U_2) $ can be viewed as
a totally real submanifold of $\C^{m+n}$. The function $G=
  F_1-  F_2$ extends holomorphically to $(U_1 \cap U_2)$ as a
function of $  w$ and is
null on  $(\U_1 \cap \U_2)$, and is thus null where it is
defined.
 
\medskip

Using the claim, we conclude that the function $S(f)$ extends
holomorphically
as a function of $  w$ to a neighborhood $\UU$ of $\Omega$ in
$\C^{m+n}$.

\medskip

\noindent {\bf (B) Holomorphy in $z$}

\bigskip

Let $  L_j \in   \L$. Since $S(f)$ extends holomorphically as a function
of $  w$, we have
\begin{equation}\label{e4.2}
  L_j(S(f))=\left [{{\partial} \over {\partial \overline {z_j}}} 
S(f) \right ]_{\Omega}.
\end{equation}

By proposition \ref{p4}, $  L_j(S(f))=S(  L_j(f))$ and since $f$ is CR and
$S(0)=0$, we have
\begin{equation}\label{e4.3}
\left [{{\partial} \over {\partial \overline {z_j}}}  
S(f) \right ]_{\Omega}=0.
\end{equation}

Note then that ${{\partial} \over {\partial \overline {z_j}}} 
S(f)$ is a holomorphic function of $  w$, for $z$ fixed. 
It vanishes on a totally real generic submanifold of $\C^{n+m}$ and thus
it is null where it is defined.

\medskip
 
 We note that $\ S(f)$ is
holomorphic on a set of the form $\B_{3\epsilon} \times \left [\R_+^{n+m} \times i
\Gamma \right ]$
where $R_+^{n+m}=\R^{n+m} \cap \{t_1>0\}$ and $\Gamma \subset
\R^{n+m}$,
thus by theorem 2.5.10 in \cite{[Ho]}, we conclude that if $ch\Gamma$
stands
for the convex hull of $\Gamma$ then  $S(f)$ is
holomorphic on  $\B_{3\epsilon} \times \left [\R_+^{n+m} \times i ch\Gamma \right
] $.

\bigskip

\bigskip

To conclude the proof of proposition \ref{Prop1} we need to show that
 the holomorphic extension obtained has slow growth and thus its
 boundary
value agrees with
our original CR distribution $f$. $S(f)$ is holomorphic in some wedge $\W$
and admits a boundary value when approaching $\n \times \{0\}$ along
$\n \times \R^{n} \cap \{t_1>0\}$, hence $S(f)$ has at most slow
growth when approaching $\n \times \{0\}$ along
$\n \times \R^{n} \cap \{t_1>0\}$. By continuity, it has at most slow
 growth
on any wedge $\W^0 \subset \W$, hence
 $S(f)$ admits a
boundary value on $\n \times \{0\}$, which by uniqueness of boundary
 value
is $f$.
This concludes the proof of proposition \ref{Prop1}. $\blacksquare$

\section{Proof of the Main Results and Remarks}

\subsection {Proof of the Results on Decomposable CR Distributions}

\noindent {\bf Proof of Theorem \ref{t1}.} Recall that 
$  N=\{(  \n,h(  \n))\}$.
Consider the CR map $h: 
\n \to \C^n$. By theorem \ref{Th2}, each $h_j$ extends holomorphically
to some wedge $\W_j$. Set $\W =\cap_{j=1}^n \W_j$.
Define $F:(  \n,0) \to (\n,\kappa h(  \n))$ where $\kappa \in \R^*$ by
$$F(z,w',w'')=(z,w',w''+\kappa h(z,w')).$$
Clearly, there exists $\kappa\not = 0$ so that on  $\overline \W$, the
Jacobian
of $F$ is non zero, without loss of generality, we can assume that
$\kappa=1$.
 Hence $F$ is a biholomorphism from $\W$
to
$F(\W)$ extending smoothly to a diffeomorphism from
$\overline \W$ to $F(\overline \W)$. Since $dh(0)=0$, $F$ is tangent
to the identity at the origin, hence there exists $\W'$ a wedge of
direction $v$ such that $\W'\subset F(\W)$. We then conclude that
 any decomposable CR distribution on $N$ extends
holomorphically to a complex transversal wedge of direction $v$.
 
To prove (B), we note that
$f_j=w''_j-h_j$ are holomorphic in  $\W$ and, since $dh_j(0)=0$, 
we conclude that
$$d(w''_1-h_1)\wedge...\wedge d(w''_n-h_n)\not =0~~{\rm on} \W.$$
Each $f_j$ vanishes to order one on $N$ since
 $F^{-1}(z,w',w'')=(z,w',w''-h)=(z,w',t)$. $\blacksquare$

\bigskip

\noindent {\bf Proof of Corollary \ref{t3}.} Let $M$ and $N$ be as in
the hypothesis of theorem \ref{t3}. After a linear of variables, we
may assume that $p_0=0$ and that near the origin, $M$ is parametrized
by
$$M=\{z,u+iv(z,u):(z,u) \in \C^k \times \R^{p-k}\}.$$
By the implicit function theorem, we may assume that $N$ is given
as a subset of $M$ by

$$
\begin{cases}
u_{p-k-n}=\mu_1(z,u_1,...,u_{p-k-n-1}),...,
u_{p-k}=\mu_n(z,u_1,...,u_{p-k-n-1}),\\
\mu(0)=d\mu(0)=0.
\end{cases}
$$
Denote by $s=(u_1,...,u_{p-k-n-1}) \in \R^m$ and
$t=(u_{p-k-n},...,u_{p-k}) \in \R^n$. Setting $t'=t-\mu$, in the $(z,s,t')$
coordinates,
we have $N$ given as a subset of $M$ by $t'=0$ and
$$N=\{(z,w'(z,s),h(z,w'):(z,s)\in \C^k  \times \R^m\},$$
where $h$ is a CR map from $\n:=\{z,w'(z,s)\}$. We can now apply 
 theorem
\ref{t1} to obtain the CR extension as the restriction of the
holomorphic
extension of $f$ to $\W \cap M$.
 The second part follows in the same manner.
$\blacksquare$

\subsection {The non Decomposable Case}

\noindent {\bf Tr\'epreau's Example}

\bigskip

We wish to note that on the particular class of manifolds of CR
dimension $k$ containing through the origin $\C^k$, we do get a complex
transversal extension result for any CR functions. Tr\'epreau's
original
example of a CR submanifold of $\C^3$ which admits non decomposable CR
functions, is the following
\begin{equation}\label{e111}
\m=\{(z,s_1+is_2|z|^2,s_2-is_1|z|^2):(z,s_1,s_2)\in \C \times \R
\times \R\} 
\end{equation}
On this particular example, we do get holomorphic extension to {\bf complex
transversal} wedges of {\bf any} CR distribution.

\begin{proposition}\label{lader} Let $M=\{(\m,h(\m)\}$ be a non
generic CR submanifold of $\C^L$ where $\m$ is given by \ref{e111},
then the conclusions of theorem \ref{t1}
hold for  CR distributions.
\end{proposition}

\noindent{\bf Proof of Proposition \ref{lader}.} As seen previously, there
is no loss of generality in assuming that 
$M=\{(\m,0)\} \subset \C^{k+m}\times \C^n$ and that if $v$ is a
complex transversal vector to $M$ at the origin, then $v=(1,0,....,0)
\subset \R^n \times \{0\} \subset \C^n$. Let $f$ be a $\CC^0$ CR
function on $M$, consider $\m\times
\R^n$ and the associated elliptic operator $\Delta_{\m\times
\R^n}$. As previously, we solve a Dirichlet problem with boundary data
$f$. Denote the solution $S(f)$, by our previous arguments $S(f)$ is
holomorphic in the variables $w$ on a neighborhood of $\Omega$ denoted
$\V(\Omega)$. Note that since $\m$ contains through the origin
  $\C\times \{0\}$ then $\V(\Omega)$  contains 
$\left (\C\times \{0\}\times \R^n\right )\cap \{t_1>0\}$. Let
$\U$ be a neighborhood of the origin in $\C_z$ and let $\xi$ be a
smooth cut off function which is one on $\U$ and zero on a
neighborhood
of $\overline \U$. 
Set $F=\xi \overline \partial S(f)$ a $(0,1)$
form on a $\V(\Omega)$ with continuous coefficient. Since $S(f)$ is
holomorphic in $w$, we have $\overline \partial S(f)={{\partial}\over
{\partial {\overline z}}}S(f)$. Since $f$ is CR,  we get 
$F|_{(\m \cap \U) \times\R^n}=0$.  Let
 $u={{\partial}\over
{\partial {\overline z}}}F\!*\!{{1}\over{\pi z}}
$ then $u$ is holomorphic in
$w$
and $u|_{\m\times
\R^n}=0$. The desired extension for $f$ is thus obtained by
considering
$S(f)-u$. $\blacksquare$

\bigskip

We can generalize this result to any CR submanifold of CR dimension
$k$ which contains through a point $p_0$ $\C^k\times \{0\}$. We state
without proving the next result, since the proof of this result is
identical to the proof of  proposition \ref{lader} one has to replace
${{1}\over{\pi z}}$ by an appropriate integral kernel solving
$\overline \partial_z$
on a ball in $\C^k$.

\begin{theorem}\label{der} Let $M$ be a non generic smooth ($\CC^{\infty}$) CR
submanifold of $\C^L$ of CR dimension $k$. If the reunion of the
 CR orbits through 
$ p_0\in M$ is a complex manifold
 then the conclusions of theorem \ref{t1}
 hold for  CR distributions.
\end{theorem}

We thus get as a corollary

\begin{corollary} (Theorem \ref{th}) Let $N=\{(\n,h(\n))\} $ 
be a smooth ($\CC^{\infty})$ non generic CR submanifold of $\C^{k+m+n}$ 
such that $\n \subset \C^{k+m}$ is a hypersurface. If $f$ is a CR
distribution on $N$ then for any $p \in N$ and any $v$ complex
transversal to $N$ at $p$, there exists a wedge
$\W$ of direction $v$ whose edge contains a neighborhood of 
$p$ in $N$
and $F\in \O(\W)$ such that
$F|_N=f$.
\end{corollary}

\bigskip 

\noindent {\bf Proof of the corollary} To prove the corollary, we need
the following result.

\begin{theorem}\label{th1} Let $\n$ be a generic submanifold of $\C^n$
minimal
as some $p\in \n$. Then any CR distribution near $p$ admits a unique
holomorphic extension to some wedge $\W$ with edge $\n$.
\end{theorem}
\noindent {\bf Proof of theorem \ref{th1}.} This result is well known
so we will only sketch the proof. Set $k=CRdim (\n)$, choose
coordinates in $\C^n$ so that $\C^n=\C^k_z \times \C^{n-k}_w$.
 Construct the vector fields $R_j$ as in lemma \ref{l3}.
 Denote by $\{L_j\}_{j=1}^{k}$ a basis for the CR vector fields.
By proposition \ref{p4} (using the terminology of Baouendi
and Treves) the system of vector fields $R_j$ is normal at $p$
(definition 4.1 in \cite{[Ba-Tr]}). Let $f$ be a CR distribution on
$\n$, by theorem 4.1 in
\cite{[Ba-Tr]} there exists $u\in \CC^0$ and $\ell \in \NN$ such that
$$\left [\Delta_{\n \times \R^n}\right ]^{\ell} u=f,$$
$$L_j(u)=0.$$
By Tumanov's theorem
 $u$ extends holomorphically to some wedge $\W$, denote by $U$
the holomorphic extension of $u$ to $\W$. The
holomorphic
extension of $f$ is then given by 
$(\sum_{j=1}^{n-k}{{\partial^2} \over {\partial w_j^2}})^{\ell}U$. $\blacksquare$

\bigskip

We now proceed with the proof of the corollary.
Let $p\in N$, $p=(p',h(p'))$,
if $\n$ is minimal at $p'$, then any CR distribution is 
decomposable (it extends holomorphically to a wedge by the above
theorem), hence
by theorem \ref{t1} we are done. If $\n$ is not minimal at $p'$, then
it contains a proper submanifold $\n'$ of same CR dimension, since
$\n$
is a hypersurface, $\n'$ is a complex manifold, hence by theorem
\ref{der}
we are done. $\blacksquare$

\bigskip

We can now construct functions holomorphic in a complex transversal
wedge
vanishing to infinite order on a non generic CR submanifold of $\C^L$.

\begin{corollary}\label{coro} Let $N=\{(\n,h(\n))\}$ be a non generic
 smooth CR submanifold of $\C^L$ such that either $\n$ is a
 hypersurface
or $h$ is decomposable at $p'\in \n$. Then for any $v$ complex
 transversal
to $N$ at $p=(p',h(p'))$ there exists a wedge $\W$ of direction $v$
and $g \in \O(\W)$ such that $g \not \equiv 0$ and $g$ vanishes to
 infinite
order on $N$.
\end{corollary}
\noindent {\bf Proof of Corollary \ref{coro}.} Choose coordinates
$(z,w',w'') \in \C_z^k\times \C_{w'}^m\times \C_{w''}^n$ such that
$\n \subset  \C_z^k\times \C_{w'}^m$. After a linear change of
variables, we can assume that $v$, the complex transversal vector is
given by $v=(0,0,v'') \in \C^{k}\times \C^m \times \C^n$ with
$v''=(1,0,...,0)$.  Consider the function $\mu$ defined by
$$\mu=
\begin{cases}
e^{-{{1}\over{w''_1}}}~~~~w''_1\not=0,\\
0~~~~w''_1\not=0.
\end{cases}
$$
Then $\mu$ is holomorphic on any wedge with direction $v$ off edge
$\n\times\{0\}$ and vanishes
to infinite order on $\n\times\{0\}$. Consider $F$ the biholomorphism
constructed in the proof of theorem \ref{t1},
$F(z,w',w'')=(z,w',w''-h(z,w'))$.
The desired function $g$ is then given by $g=\mu(F^{-1}(z,w',w''))$.
$\blacksquare$

\bigskip

\bigskip

{\bf Note.} The function $g$ does not extend holomorphically to any
neighborhood of $N$. This situation differs greatly from the
holomorphic
extension situation in the generic case. Indeed, let $M$ be a generic
manifold where any CR function in a neighborhood of some point $p\in
M$ extends holomorphically to a neighborhood of $p$. Let $\W$ be a
wedge attached to $M$, then any function which is holomorphic in $\W$
and
at least $\CC^1$ on $M$ extends holomorphically to a full neighborhood
of $p$. However, if we now consider $\M=\{(M,h(M))\}$ then the
function
$g$ in corollary \ref{coro} does not extend holomorphically to a full
neighborhood of $\M$.

\bigskip

\bigskip

\begin{remark}

\noindent  In \cite{[Ei]} we constructed an example of an abstract CR structure
where there is no CR extension. It is even easier to construct an
example where there is no non trivial CR function vanishing on $N$.
\end{remark}
\noindent {\bf Example.} 
Let $L$ be a real analytic  vector field (the Lewy operator
for example)  of the form
$$L={{\partial} \over {\partial \bar z}}+f(z,s){{\partial} \over
{\partial s}},$$
that is not solvable (H\"ormander's theorem \cite{[Ho]}). Let $g$ be a
smooth
function not in the image of $L$, then define the abstract CR
structure
$(M,\L)$ where $\L$ is given by
$$\L=L+tg{{\partial} \over
{\partial t}},$$
then the equation $\L(tu)=0$ is not solvable. Indeed, suppose it was,
then
we would have 
$$t\L(u)+u\L(t)=t(\L(u)+ug)=0,$$
decomposing $u=u_0+tu_1$ where $u_0=u_0(z,s)$ we obtain
$$L(Log(u_0))+g=0.$$

By real analycity of $L$, we see that there are plenty of non trivial
functions in the kernel of $\L$.

\bigskip

\noindent Nicolas Eisen\\
D\'epartement de Math\'ematiques, UMR 6086 CNRS\\
Universit\'e de Poitiers\\
eisen@math.univ-poitiers.fr

\end{document}

%% file: p2.pstex_t
\begin{picture}(0,0)%
\includegraphics{p2.pstex}%
\end{picture}%
\setlength{\unitlength}{1906sp}%
\begingroup\makeatletter\ifx\SetFigFont\undefined
\def\x#1#2#3#4#5#6#7\relax{\def\x{#1#2#3#4#5#6}}%
\expandafter\x\fmtname xxxxxx\relax \def\y{splain}%
\ifx\x\y   
\gdef\SetFigFont#1#2#3{%
  \ifnum #1<17\tiny\else \ifnum #1<20\small\else
  \ifnum #1<24\normalsize\else \ifnum #1<29\large\else
  \ifnum #1<34\Large\else \ifnum #1<41\LARGE\else
     \huge\fi\fi\fi\fi\fi\fi
  \csname #3\endcsname}%
\else
\gdef\SetFigFont#1#2#3{\begingroup
  \count@#1\relax \ifnum 25<\count@\count@25\fi
  \def\x{\endgroup\@setsize\SetFigFont{#2pt}}%
  \expandafter\x
    \csname \romannumeral\the\count@ pt\expandafter\endcsname
    \csname @\romannumeral\the\count@ pt\endcsname
  \csname #3\endcsname}%
\fi
\fi\endgroup
\begin{picture}(9462,6130)(406,-6014)
\put(5131,-5956){\makebox(0,0)[lb]{\smash{\SetFigFont{6}{7.2}{rm}{\color[rgb]{0,0,0}$\B_{3\epsilon} \times B_{2\epsilon}$}%
}}}
\put(406,-5461){\makebox(0,0)[lb]{\smash{\SetFigFont{6}{7.2}{rm}{\color[rgb]{0,0,0}$\B_{3\epsilon} \times B_{\epsilon}$}%
}}}
\put(7651,-2446){\makebox(0,0)[lb]{\smash{\SetFigFont{6}{7.2}{rm}{\color[rgb]{0,0,0}$\n$}%
}}}
\put(8056,-376){\makebox(0,0)[lb]{\smash{\SetFigFont{6}{7.2}{rm}{\color[rgb]{0,0,0}$\ti \n$}%
}}}
\put(3151,-61){\makebox(0,0)[lb]{\smash{\SetFigFont{6}{7.2}{rm}{\color[rgb]{0,0,0}$\W$}%
}}}
\put(1576,-1276){\makebox(0,0)[lb]{\smash{\SetFigFont{6}{7.2}{rm}{\color[rgb]{0,0,0}$\n$}%
}}}
\put(6256,-16){\makebox(0,0)[lb]{\smash{\SetFigFont{6}{7.2}{rm}{\color[rgb]{0,0,0}$\n=\ti \n$}%
}}}
\put(8236,-1186){\makebox(0,0)[lb]{\smash{\SetFigFont{6}{7.2}{rm}{\color[rgb]{0,0,0}$\ti \n=\C^k\times (\R^m +id)$}%
}}}
\put(8416,-5596){\makebox(0,0)[lb]{\smash{\SetFigFont{6}{7.2}{rm}{\color[rgb]{0,0,0}$\B_{3\epsilon} \times B_{3\epsilon}$}%
}}}
\end{picture}

%% file: p4.pstex_t
\begin{picture}(0,0)%
\includegraphics{p4.pstex}%
\end{picture}%
\setlength{\unitlength}{1657sp}%
\begingroup\makeatletter\ifx\SetFigFont\undefined
\def\x#1#2#3#4#5#6#7\relax{\def\x{#1#2#3#4#5#6}}%
\expandafter\x\fmtname xxxxxx\relax \def\y{splain}%
\ifx\x\y   
\gdef\SetFigFont#1#2#3{%
  \ifnum #1<17\tiny\else \ifnum #1<20\small\else
  \ifnum #1<24\normalsize\else \ifnum #1<29\large\else
  \ifnum #1<34\Large\else \ifnum #1<41\LARGE\else
     \huge\fi\fi\fi\fi\fi\fi
  \csname #3\endcsname}%
\else
\gdef\SetFigFont#1#2#3{\begingroup
  \count@#1\relax \ifnum 25<\count@\count@25\fi
  \def\x{\endgroup\@setsize\SetFigFont{#2pt}}%
  \expandafter\x
    \csname \romannumeral\the\count@ pt\expandafter\endcsname
    \csname @\romannumeral\the\count@ pt\endcsname
  \csname #3\endcsname}%
\fi
\fi\endgroup
\begin{picture}(9068,3769)(485,-6239)
\put(3691,-3346){\makebox(0,0)[lb]{\smash{\SetFigFont{5}{6.0}{rm}{\color[rgb]{0,0,0}$\gamma_2$}%
}}}
\put(3286,-5416){\makebox(0,0)[lb]{\smash{\SetFigFont{5}{6.0}{rm}{\color[rgb]{0,0,0}$\gamma_1$}%
}}}
\put(1261,-5506){\makebox(0,0)[lb]{\smash{\SetFigFont{5}{6.0}{rm}{\color[rgb]{0,0,0}$B_{3\epsilon}$}%
}}}
\put(1306,-4876){\makebox(0,0)[lb]{\smash{\SetFigFont{5}{6.0}{rm}{\color[rgb]{0,0,0}$B_{2\epsilon}$}%
}}}
\put(1981,-2716){\makebox(0,0)[lb]{\smash{\SetFigFont{5}{6.0}{rm}{\color[rgb]{0,0,0}$t_1$}%
}}}
\put(5581,-6181){\makebox(0,0)[lb]{\smash{\SetFigFont{5}{6.0}{rm}{\color[rgb]{0,0,0}$\ti \n \times \R^{n-1}$}%
}}}
\put(9361,-2671){\makebox(0,0)[lb]{\smash{\SetFigFont{5}{6.0}{rm}{\color[rgb]{0,0,0}$\Omega$}%
}}}
\put(9541,-5011){\makebox(0,0)[lb]{\smash{\SetFigFont{5}{6.0}{rm}{\color[rgb]{0,0,0}$\Gamma_1$}%
}}}
\put(5671,-2626){\makebox(0,0)[lb]{\smash{\SetFigFont{5}{6.0}{rm}{\color[rgb]{0,0,0}$\Gamma_2$}%
}}}
\put(586,-2941){\makebox(0,0)[lb]{\smash{\SetFigFont{5}{6.0}{rm}{\color[rgb]{0,0,0}$\omega$}%
}}}
\end{picture}